\documentclass[10pt, twoside]{scrartcl}
\usepackage{etex}
\usepackage[a4paper, top=88pt, bottom=88pt, left=72pt, right=72pt, headsep=16pt, footskip=28pt]{geometry}
\usepackage{amsmath, amssymb, amsthm, xcolor, tabularx}
\usepackage[hyphens]{url}
\usepackage{hyperref}
\hypersetup{
	hidelinks,
}
\usepackage[nameinlink]{cleveref}
\usepackage{enumitem}
\usepackage{mleftright} \mleftright
\renewcommand{\qedsymbol}{$\blacksquare$}

\usepackage{scrlayer-scrpage, xhfill}
\automark[section]{section}
\setkomafont{pageheadfoot}{\normalcolor\sffamily}
\setkomafont{pagenumber}{\normalfont\normalsize\sffamily}
\clearpairofpagestyles
\let\oldtitle\title
\renewcommand{\title}[1]{\oldtitle{#1}\newcommand{\theshorttitle}{#1}}
\newcommand{\shorttitle}[1]{\renewcommand{\theshorttitle}{#1}}
\let\oldauthor\author
\renewcommand{\author}[1]{\oldauthor{#1}\newcommand{\theshortauthor}{#1}}
\newcommand{\shortauthor}[1]{\renewcommand{\theshortauthor}{#1}}
\cohead{\xrfill[0.525ex]{0.6pt}~\theshorttitle~\xrfill[0.525ex]{0.6pt}}
\cehead{\xrfill[0.525ex]{0.6pt}~\theshortauthor~\xrfill[0.525ex]{0.6pt}}
\usepackage{natbib}
\setcitestyle{round}
\newcommand*{\arXiv}[1]{\bgroup\color{blue}\href{https://arxiv.org/abs/#1}{arXiv:#1}\egroup}
\newcommand*{\doi}[1]{\bgroup\color{blue}\href{https://dx.doi.org/#1}{doi:#1}\egroup}
\newcommand*{\email}[1]{\bgroup\color{blue}\href{mailto:#1}{#1}\egroup}
\renewcommand*{\url}[1]{\bgroup\color{blue}\href{#1}{#1}\egroup}
\newcommand{\todo}[1]{\bgroup\color{red}\bfseries#1\egroup}

\usepackage{mathtools}

\newcommand*{\argmin}{\mathop{\textup{arg\,min}}}
\newcommand*{\Complex}{\mathbb{C}}

\newcommand*{\defeq}{\coloneqq}

\newcommand*{\Naturals}{\mathbb{N}}
\newcommand*{\Normal}{\mathcal{N}}

\newcommand*{\quark}{\setbox0\hbox{$x$}\hbox to\wd0{\hss$\cdot$\hss}}
\DeclareMathOperator{\rank}{rank}

\newcommand*{\Reals}{\mathbb{R}}
\DeclareMathOperator{\trace}{tr}

\newcommand*{\norm}[1]{\Vert #1 \Vert}

\newcommand*{\Set}[2]{\left\{ #1 \middle\vert #2 \right\}}

\theoremstyle{definition}

\numberwithin{equation}{section}
\numberwithin{figure}{section}
\numberwithin{table}{section}

\renewcommand*{\paragraph}[1]{\smallskip\noindent\textbf{\textsf{#1}}\,\,}
\newcommand*{\bs}{\boldsymbol}
\newcommand*{\Krylov}{\mathcal{K}}

\begin{document}

\title{Comments on the article\\``A Bayesian conjugate gradient method''}
\shorttitle{Comments on ``A Bayesian conjugate gradient method''}

\author{%
	T.\ J.\ Sullivan\footnote{Institute of Mathematics, Freie Universit{\"a}t Berlin, Arnimallee 6, 14195 Berlin, Germany, \email{t.j.sullivan@fu-berlin.de}} \footnote{Zuse Institute Berlin, Takustra{\ss}e 7, 14195 Berlin, Germany, \email{sullivan@zib.de}}%
}
\shortauthor{T.~J.~Sullivan}

\newcommand{\TJSsays}[1]{\bgroup\color{purple}{\text{TJS says: }#1}\egroup}

\date{\today}

\maketitle

\begin{abstract}
	\paragraph{Abstract:}
	The recent article ``A Bayesian conjugate gradient method'' by Cockayne, Oates, Ipsen, and Girolami proposes an approximately Bayesian iterative procedure for the solution of a system of linear equations, based on the conjugate gradient method, that gives a sequence of Gaussian/normal estimates for the exact solution.
The purpose of the probabilistic enrichment is that the covariance structure is intended to provide a posterior measure of uncertainty or confidence in the solution mean.
This note gives some comments on the article, poses some questions, and suggests directions for further research.


\end{abstract}

\section{Overview}

In ``A Bayesian conjugate gradient method'', Cockayne, Oates, Ipsen, and Girolami add to the recent body of work that provides probabilistic/inferential perspectives on deterministic numerical tasks and algorithms.
In the present work, the authors consider a conjugate gradient (CG) method for the solution of a finite-dimensional linear system $A \bs{x}^{\ast} = \bs{b}$ for $\bs{x}^{\ast} \in \Reals^{d}$, given an assuredly invertible symmetric matrix $A \in \Reals^{d \times d}$ and $\bs{b} \in \Reals^{d}$, with $d \in \Naturals$.

The authors derive and test an algorithm, BayesCG, that returns a sequence of normal distributions $\Normal ( \bs{x}_{m} , \Sigma_{m} )$ for $m = 0, \dots, d$, starting from a prior distribution $\Normal ( \bs{x}_{0}, \Sigma_{0} )$.
This sequence of normal distributions is defined using a recursive relationship similar to that defining the classical CG method, and indeed the BayesCG mean $\bs{x}_{m}$ coincides with the output of CG upon choosing $\Sigma_{0} \defeq A^{-1}$ ---
this choice is closely related to what the authors call the ``natural prior covariance'', $\Sigma_{0} \defeq (A^{\top} A)^{-1}$ .
The distribution $\Normal ( \bs{x}_{m} , \Sigma_{m} )$ is intended to be an expression of posterior belief about the true solution $\bs{x}^{\ast}$ to the linear system under the prior belief $\Normal ( \bs{x}_{0} , \Sigma_{0} )$ given the first $m$ BayesCG search directions $\bs{s}_{1}, \dots, \bs{s}_{m}$.
Like CG, BayesCG terminates in at most $d$ steps, at which point its mean $\bs{x}_{d}$ is the exact solution $\bs{x}^{\ast}$ and it expresses complete confidence in this belief by having $\Sigma_{d} = 0$.
The convergence and frequentist coverage properties of the algorithm are investigated in a series of numerical experiments.

The field of probabilistic perspectives on numerical tasks has been enjoying a resurgence of interest in the last few years;
see e.g.\ \citet{OatesSullivan2019} for a recent survey of both historical and newer work.
The authors' contribution is a welcome addition to the canon, showing as it does how classical methods (in this case CG;
cf.\ the treatment of Runge--Kutta methods for ODEs by \citet{Schober2014}) can be seen as point estimators of particular instances of inferential procedures.
it is particularly encouraging to see contributions coming from authors with both statistical and numerical-analytical expertise, and the possibilities for generalisation and further work are most interesting.

\section{Questions and directions for generalisation}

The article raises a number of natural questions and directions for generalisation and further investigation, which Cockayne et al.\ might use their rejoinder to address.

\paragraph{Symmetry and generalisations.}
It is interesting that, for the most part, BayesCG does not require that $A$ be symmetric, and works instead with $A \Sigma_{0} A^{\top}$.
This prompts a question for the authors to consider in their rejoinder:
How does BayesCG behave in the case that $A \in \Reals^{d \times d}$ is square but not invertible?
Could one even show that the BayesCG posterior concentrates, within at most $d$ iterations, to a distribution centred at the minimum-norm solution (in some norm on $\Reals^{d}$) of the linear system?
One motivation for this question is that, if such a result could be established, then BayesCG could likely be applied to rectangular $A \in \Reals^{c \times d}$ and produce a sequence of normally-distributed approximate solutions to the minimum-norm least-squares problem, i.e.\ a probabilistically-motivated theory of Moore--Penrose pseudo-inverses.
(Recall that the Moore--Penrose pseudo-inverse $A^{\dagger} \in \Reals^{d \times c}$ can be characterised as the solution operator for the minimum-norm least squares problem, i.e.,
\[
	A^{\dagger} \bs{b} = \argmin \Set{ \norm{ \bs{x} }_{\Reals^{d}} }{ \bs{x} \in \argmin \Set{ \norm{ A \bs{x}' - \bs{b} }_{\Reals^{c}} }{ \bs{x}' \in \Reals^{d} } }
\]
for the usual norms on $\Reals^{c}$ and $\Reals^{d}$.)

The authors could also comment on whether they expect BayesCG to generalise easily to infinite-dimensional Hilbert spaces, analogously to CG \citep{Fortuna1977, Fortuna1979, Malek2015}, since experience has shown that analysis of infinite-dimensional statistical algorithms can yield powerful dimension-independent algorithms for the finite-dimensional setting \citep{Cotter2013, Chen2018}.
A more esoteric direction for generalisation would be to consider fields other than $\Reals$.
The generalisation of BayesCG to $\Complex$ would appear to be straightforward via a complex normal prior, but how about fields of finite characteristic?

\paragraph{Conditioning relations.}
Could the authors use their rejoinder to provide some additional clarity about the relationship of BayesCG to exact conditioning of the prior normal distribution $\Normal (\bs{x}_{0}, \Sigma_{0})$, and in particular whether BayesCG is exactly replicating the ideal conditioning step, approximating it, or doing something else entirely?

To make this question more precise, fix a weighted inner product on $\Reals^{d}$, e.g.\ the Euclidean, $A$-weighted, $\Sigma_{0}$-weighted, or $A \Sigma_{0} A^{\top}$-weighted inner product.
With respect to this inner product, let $P_{m}$ be the (self-adjoint) orthogonal projection onto the Krylov space $\Krylov_{m}$ and $P_{m}^{\perp} \defeq I - P_{m}$ the (self-adjoint) orthogonal projection onto its orthogonal complement.
The product measure
\[
	\mu_{m} \defeq \Normal(P_{m} \bs{x}^{\ast}, 0) \otimes \Normal ( P_{m}^{\perp} \bs{x}_{0} , P_{m}^{\perp} \Sigma_{0} P_{m}^{\perp} )
\]
on $\Krylov_{m} \oplus \Krylov_{m}^{\perp}$ is also a normal distribution on $\Reals^{d}$ that expresses complete confidence about the true solution to the linear system in the directions of the Krylov subspace and reverts to the prior in the complementary directions;
$\mu_{0}$ is the prior, and $\mu_{d}$ is a Dirac on the truth.
The obvious question is, for some choice of weighted inner product, does BayesCG basically output $\mu_{m}$, but in a clever way?
Or is the output of BayesCG something else?

\paragraph{Uncertainty quantification: cost, quality, and terms of reference.}
Cockayne et al.\ point out that the computational cost of BayesCG is a small multiple (a factor of three) of the cost of CG, and this would indeed be a moderate price to pay for high-quality uncertainty quantification.
However, based on the results presented in the paper, the UQ appears to be quite poorly calibrated.
One could certainly try to overcome this shortcoming by improving the UQ.
However, an alternative approach would be to choose a prior covariance structure $\Sigma_{0}$ that aggressively sparsifies and hence accelerates the linear-algebraic operations that BayesCG performs (particularly lines 7, 9, and 11 of Algorithm~1), while preserving the relatively poor UQ.
Does this appear practical, in their view?

In fact, the whole question of ``the UQ being well calibrated'' is somewhat ill defined.
Some might argue that, since $\bs{x}^{\ast}$ is deterministic, there is no scope for uncertainty in this setting.
However, It certainly makes sense to ask --- as the authors do in Section~6.1 -- whether, when the problem setup is randomised, the \emph{frequentist} coverage of BayesCG lines up with that implied by the randomisation of the problem.
The statistic $Z$ that Cockayne et al.\ introduce is intersting, and already captures several scenarios in which the UQ is well calibrated and several in which it is not;
I strongly encourage Cockayne et al.\ to address in their rejoinder the \emph{Bayesian} accuracy of BayesCG, e.g., to exhaustively sample the true posterior on $\bs{x}^{\ast}$ given $\bs{y}_{1}, \dots, \bs{y}_{m}$ and see whether this empirical distribution is well approximated by the BayesCG distribution $\Normal ( \bs{x}_{m}, \Sigma_{m} )$.

As a related minor question, can BayesCG be seen as an (approximate Gaussian) \emph{filtering} scheme for the given prior $\Normal(\bs{x}_{0}, \Sigma_{0})$ and the filtration associated to the data stream $\bs{y}_{1}, \bs{y}_{2}, \dots$?

\paragraph{Precision formulation.}
As formulated, BayesCG expresses a Gaussian belief about the solution of the linear system in terms of mean and covariance;
Gaussian measures can also be expressed in terms of their mean and precision.
Do the authors have a sense of whether the BayesCG algorithm be formulated as a sequence of precision updates, or does the singularity of the covariance in the Krylov directions essentially forbid this?

One motivation for seeking a precision formulation would be to render the ``natural prior'' $\Sigma_{0} \defeq (A^{\top} A)^{-1}$ of Section 4.1 more tractable, since this prior has an easily-accessible precision while accessing its covariance involves solving the original linear problem.

As the authors note, their ``natural prior'' is closely related to one introduced by \citet{Owhadi2015} and applied in \citet{Cockayne2016}.
It seems that working with images of Gaussian white noise, such as this ``natural prior'', is presently producing considerable analytical and computational strides forward \citep{Owhadi2017, Chen2018}, and so this seems to be a topic worth further attention in the statistical community as a whole.

\section{Minor comments}

\paragraph{Rank and trace estimates.}
Proposition~3, which states that $\trace ( \Sigma_{m} \Sigma_{0}^{-1} ) = d - m$, seems to miss the point slightly.
It would be good to have a companion results to the effect that $\rank \Sigma_{m} = d - m$, and a more quantitative result for the trace such as $\trace \Sigma_{m} \leq C_{\Sigma_{0}} ( d - m )$ for some constant $C_{\Sigma_{0}} \geq 0$ depending only on $\Sigma_{0}$.

\paragraph{Posterior and Krylov spaces.}
It seems natural to ask whether the posterior nullspace $\ker \Sigma_{m}$ coincides with the Krylov space $\Krylov_{m}$.
Put another way, is the posterior column space the same as the orthogonal complement of the Krylov space, in the Euclidean or $A$-weighted inner product?

\paragraph{Square roots.}
Is the square root $M^{1/2}$ introduced just before Proposition~5 required to be unique?
Does it even matter whether a (unique) symmetric positive semidefinite square root or a Cholesky factor is chosen?
This is related to the above discussion of symmetry and generalisations.

\paragraph{Interpretation of termination criteria.}
In Section 5, the authors refer to ``probabilistic termination criteria''.
As a matter of semantics, the termination criterion that they propose is in fact deterministic, albeit based on a probabilistic interpretation of the algorithmic quantities.


\section*{Acknowledgements}
\addcontentsline{toc}{section}{Acknowledgements}
TJS is supported by the Freie Universit\"at Berlin within the Excellence Strategy of the German Research Foundation (DFG), including the ECMath/MATH+ transition project CH-15 and project TrU-2 of the Berlin Mathematics Research Center MATH+ (EXC-2046/1, project ID 390685689).

\bibliographystyle{abbrvnat}
\addcontentsline{toc}{section}{References}
\bibliography{references.bib}

\end{document}